\documentclass[12pt,a4paper,reqno]{amsart}
\usepackage[english]{babel}
\usepackage{amssymb,latexsym,amsfonts,amsthm,upref,amsmath}
\usepackage[foot]{amsaddr}
\usepackage[margin=1in]{geometry} 
\usepackage[dvipsnames,x11names]{xcolor}
 \definecolor{myblue}{HTML}{003399}
\usepackage{hyperref}
\hypersetup{colorlinks,citecolor=myblue,filecolor=black,linkcolor=myblue,urlcolor=myblue}
\usepackage{enumerate}  
\usepackage{tikz}
\usepackage{float}
\usepackage{cleveref}
\usepackage{mathtools}
\usepackage{cite}
\usepackage{filecontents}
\makeatletter
\newcommand{\leqnomode}{\tagsleft@true}
\newcommand{\reqnomode}{\tagsleft@false}
\newcommand{\cev}[1]{\reflectbox{\ensuremath{\vec{\reflectbox{\ensuremath{#1}}}}}}
\makeatother
\newtheorem*{thm*}{Theorem}
\newtheorem*{lem*}{Lemma}
\newtheoremstyle{prim}{}{}{\normalfont}{}{\bfseries}{.}{ }{}
\newtheoremstyle{stil}{}{}{\slshape}{}{\bfseries}{.}{ }{}
\theoremstyle{stil}
\newtheorem{thm}{Theorem}[section]
\newtheoremstyle{defi}{}{}{}{}{\bfseries}{.}{ }{}
\theoremstyle{defi}
\newtheorem{defn}[thm]{Definition}
\theoremstyle{defi}
\newtheorem{rem}[thm]{Remark}
\theoremstyle{stil}
\newtheorem*{mthm*}{Main Theorem}
\newtheorem*{kor*}{Corollary}
\newtheorem{pro}[thm]{Proposition}
\theoremstyle{stil}
\newtheorem{lem}[thm]{Lemma}
\theoremstyle{stil}
\newtheorem{kor}[thm]{Corollary}
\theoremstyle{prim}
\newtheorem{ex}[thm]{Example}
\newenvironment{prf}{\noindent \textit{Proof.}}{\null\hfill$\qed$\hskip
2mm\vskip 2mm}
\newenvironment{prfn}{\noindent \textit{Proof.}} 

\newcommand{\Yo}{ {\rm Y}_{h}^{\rm tw}(\mathfrak{o}_{N})}
\newcommand{\Xo}{ {\mathcal{Y}}_{h}^{\rm tw}(\mathfrak{o}_{N})}

\newcommand{\CC}{\mathbb{C}}
\newcommand{\Sc}{\mathcal{S}}

\newcommand{\Vcgl}{\mathcal{V}^c(\mathfrak{gl}_N)}

\newcommand{\R}{\wvr{R}}

\newcommand{\RR}{\wnr{\wvr{R}}}
\newcommand{\RRR}{\underline{R}}
\newcommand{\ZZ}{\mathbb{Z}}
\newcommand{\wvr}{\overline}

\newcommand{\wnr}{\underline}
\newcommand{\vac}{\mathop{\mathrm{\boldsymbol{1}}}}

\newcommand{\tr}{\mathop{\mathrm{tr}}}

\newcommand{\ot}{\otimes}
\newcommand{\ts}{\hspace{1pt}}
\newcommand{\qdet}{ \mathop{\rm qdet} }
\newcommand{\sgn}{ \mathop{\rm sgn}}

\newcommand{\ndo}{\mathop{\mathrm{End}}}
\newcommand{\om}{\mathop{\mathrm{Hom}}}

\newcommand{\diag}{\mathop{\mathrm{diag}}}
\newcommand{\cdotrl}{\mathop{\hspace{-2pt}\underset{\text{RL}}{\cdot}\hspace{-2pt}}}
\newcommand{\cdotlr}{\mathop{\hspace{-2pt}\underset{\text{LR}}{\cdot}\hspace{-2pt}}}

\newcommand{\iotan}{\mathop{\iota_{z_1,\ldots ,z_n}}}
\newcommand{\iotadva}{\mathop{\iota_{z_1, z_2}}}
\newcommand{\iotau}{\mathop{\iota_{u}}}
\newcommand{\iotaze}{\mathop{\iota_{z}}}

\newcommand{\iotazepm}{\mathop{\iota_{ z, x^{\pm 1}}}}
\newcommand{\iotazepmrev}{\mathop{\iota_{x^{\pm 1},z}}}
\newcommand{\iotaop}{\mathop{\iota}}
\newcommand{\iotaopdvanula}{\mathop{\iota_{z_2,z_0}}}
\newcommand{\iotasigma}{\mathop{\iota_{u_{\sigma_1} ,\ldots ,u_{\sigma_n}}}}

\newcommand{\Ur}{{\rm{U}}(R)}
\newcommand{\wht}{\widehat}

\newcommand{\fand}{\quad\text{and}\quad}
\newcommand{\Fand}{\qquad\text{and}\qquad}

\newcommand{\non}{\nonumber}
\newcommand{\beq}{\begin{equation}}
\newcommand{\eeq}{\end{equation}}
\newcommand{\ben}{\begin{equation*}}
\newcommand{\een}{\end{equation*}}

\makeatletter
\def\smalloverbrace#1{\mathop{\vbox{\m@th\ialign{##\crcr\noalign{\kern3\p@}%
  \tiny\downbracefill\crcr\noalign{\kern3\p@\nointerlineskip}%
  $\hfil\displaystyle{#1}\hfil$\crcr}}}\limits}
\makeatother

\makeatletter
\def\smallunderbrace#1{\mathop{\vtop{\m@th\ialign{##\crcr
   $\hfil\displaystyle{#1}\hfil$\crcr
   \noalign{\kern3\p@\nointerlineskip}%
   \tiny\upbracefill\crcr\noalign{\kern3\p@}}}}\limits}
\makeatother

\begin{document}

\title{Associating deformed  \texorpdfstring{$\phi$}{phi}-coordinated   modules for the quantum affine vertex algebra  with orthogonal twisted $h$-Yangians }

\author{Lucia Bagnoli}
\author{Slaven Ko\v{z}i\'{c}}
\address[L. Bagnoli and S. Ko\v{z}i\'{c}]{Department of Mathematics, Faculty of Science, University of Zagreb,  Bijeni\v{c}ka cesta 30, 10\,000 Zagreb, Croatia}
\email{lucia.bagnoli@math.hr}
\email{kslaven@math.hr}

\begin{abstract}
We consider the Etingof--Kazhdan quantum vertex algebra $\mathcal{V}^c(\mathfrak{gl}_N)$ associated with the trigonometric $R$-matrix of type $A$. 
By combining   Li's theory of $\phi$-coordinated modules and the
ideas from our previous paper, we introduce the notion of deformed $\phi$-coordinated quantum vertex algebra module.
We show that the orthogonal twisted $h$-Yangians and restricted modules for the generalized orthogonal twisted $h$-Yangians  can be equipped with the structure of (truncated)
deformed  $\phi$-coordinated  $\mathcal{V}^c(\mathfrak{gl}_N)$-module 
and   demonstrate its applications. 
\end{abstract}

 \maketitle

\allowdisplaybreaks

 \section{Introduction}
 \numberwithin{equation}{section}

The definition of quantum vertex algebra was given by Etingof and Kazhdan \cite{EK5}. Moreover, they   constructed quantum affine vertex algebras associated with rational, trigonometric and elliptic $R$-matrices. Another major step in the development  of general quantum vertex algebra theory was Li's introduction of $\phi$-coordinated modules \cite{Liphi}. In particular, it enabled the association of quantum vertex algebras with quantum affine algebras and their generalizations;  for more information see the papers by  Jing,   Kong,  Li, Tan \cite{JKLT} and Kong \cite{Kong} and references therein.

In this paper, we continue the research of Molev and the second author \cite{KM,K} on the Etingof--Kazhdan quantum affine vertex algebra $\Vcgl$ associated with  the trigonometric $R$-matrix of type $A$. Our  goal is to extend the notion of {\em deformed module}, from our previous paper \cite{BK}, to the   $\phi$-coordinated module setting. The definition of deformed module is characterized by more general forms of the weak associativity and the $\mathcal{S}$-locality of the module map $Y_W(\cdot ,z)$. They are governed by a certain map $\nu=\nu(z,x)$, defined on the tensor square $V\ot V$ of the corresponding quantum vertex algebra $V$, so that they possess the form\footnote{For a more precise formulation, see \cite[Def. 3.5]{BK}.}
\beq\label{intro_assoc}
 Y_W(u,z_0 +z_2)\ts Y_W(v,z_2)\ts w 
\sim    Y_W\big(Y ( z_0 )\ts \nu( z_0,z_2)(u\ot v),z_2\big)\ts w  
\eeq
and
\beq\label{intro_loc}
Y_W(z_1)\big(1\otimes Y_W(z_2)\big) \mu( z_1-z_2 ,z_2)(u\otimes v) \ot w
\sim Y_W(v,z_2)  Y_W(u,z_1)w ,
\eeq
respectively, where $Y(\cdot ,z)$ is the vertex operator map of $V$, $\mathcal{S}$   the braiding of $V$, $W$ a deformed $V$-module, $u,v\in V$, $w\in W$ and
\beq\label{intro_mu}
\mu(z,x) = \nu(z,x)^{-1}\Sc(z)\ts \nu_{21}(-z,x+z) .
\eeq
We refer to a pair $(\mu,\nu)$  of such maps as a {\em compatible pair}. Clearly, by setting $\nu=1$, the identity in \eqref{intro_mu} yields $\mu=\Sc$, so the properties \eqref{intro_assoc} and \eqref{intro_loc} turn to the ordinary  weak associativity and $\mathcal{S}$-locality for the quantum vertex algebra module map, as given by Li \cite{Li}.

The paper is organized as follows. In Section \ref{section_01}, we recall the Etingof--Kazhdan construction of the quantum affine vertex algebra $\Vcgl$  associated with the trigonometric $R$-matrix of type $A$.
In Section \ref{section_02}, we find the suitable $\phi$-coordinated  analogue of the  notion of compatible pair, which is in tune with the associate $\phi(z_0,z_2)=z_2e^{z_0}$; we refer to it as a {\em multiplicative compatible pair}\footnote{It is worth noting that, as implicitly suggested by the preceding paragraph, the   compatible pairs from \cite{BK} are in tune with the trivial associate $\phi(z_0,z_2)=z_0 + z_2$.}. In general, $\phi$-coordinated modules can be regarded for any associate $\phi(z_0,z_2)\in\CC((z_0))[[z_2]]$ of  a one-dimensional additive formal group; see \cite{Liphi} for details. However, we chose this particular series $\phi$ due to its close connection with   quantum affine algebras and, consequently, twisted $h$-Yangians, as they are of central interest in this paper. Finally, we construct an example of a 
multiplicative compatible pair $(\mu,\nu)$ over the quantum  vertex algebra $\Vcgl$.

In Sections \ref{section_03} and \ref{section_04}, we associate to a  multiplicative compatible pair $(\mu,\nu)$  the notion of {\em (truncated) $(\mu,\nu)$-deformed $\phi$-coordinated module}.
Next, we use the   pair $(\mu,\nu)$, constructed in   Section \ref{section_02}, to establish a connection between the quantum vertex algebra $\Vcgl$ and the {\em orthogonal twisted $h$-Yangian} $\Yo$. 
This coideal subalgebra of the quantum loop algebra of $\mathfrak{gl}_N$ goes back to
   Molev,   Ragoucy and Sorba \cite{MRS}, where its
	presentation  in terms of generators and defining relations in
the form of reflection  equation was given; see also the recent paper \cite{Lu} by Lu for its  Gauss decomposition. 
In Section \ref{section_03}, we show that the truncated $(\mu,\nu)$-deformed $\phi$-coordinated $\Vcgl$-modules are essentially the $\Yo$-modules. Furthermore, we find a new,  vertex algebraic proof of centrality of the coefficients of  {\em Sklyanin determinant}, a quantum determinant analogue  for $\Yo$, which was originally proved in  \cite{MRS}. 
It employs, in  particular, the explicit description of the center of $\Vcgl$ at the noncritical level $c\neq -N$, which goes back to  \cite{KM}.
In Section  \ref{section_04}, we introduce the 
{\em generalized orthogonal twisted $h$-Yangian} $\Xo$ and we show that its  (suitably defined)  {\em restricted modules} are naturally equipped with the structure of $(\mu,\nu)$-deformed $\phi$-coordinated $\Vcgl$-module.

\section{Preliminaries}\label{section_01}

\subsection{Trigonometric \texorpdfstring{$R$}{R}-matrix}\label{subsection_11}
Let $N\geqslant 2$ be an integer and $h$ a formal parameter. Denote by $\CC[[h]]$ a commutative ring of formal Taylor series in $h$ with complex coefficients. Consider the  trigonometric $R$-matrix $\R(x)\in \ndo\CC^N\ot\ndo\CC^N[[h]][x]$ of type $A$,
\begin{align}
\R(x) =&\left(1-e^{-h}x\right)\sum_{i=1}^N e_{ii}\ot e_{ii} 
+ e^{-h/2} \left(1-x\right)\sum_{\substack{i,j=1\\i\neq j}}^N e_{ii}\ot e_{jj}\non\\
&+ \left(1-e^{-h}\right)x \sum_{\substack{i,j=1\\i> j}}^N e_{ij}\ot e_{ji}
+\left( 1-e^{-h}\right)\sum_{\substack{i,j=1\\i< j}}^N e_{ij}\ot e_{ji}
,\label{Rbar}
\end{align}
where $e_{ij}\in\ndo \CC^N$ are the matrix units and $e^{ah}=\sum_{k\geqslant 0} (ah)^k /k!\in\CC[[h]]$ for $a\in\CC$.
 It satisfies the {\em quantum Yang--Baxter equation}
\beq\label{YBE}
\R_{12}(x/y)\R_{13}(x)\R_{23}(y)=\R_{23}(y)\R_{13}(x)\R_{12}(x/y)  
\eeq 
and the following identity, which resembles the unitarity property,
\beq\label{unitrig}
\R_{12}(x)\ts \R_{21}(1/x) =\left(1-e^{-h}x\right)\left(1-e^{-h}x^{-1}\right).
\eeq 
Throughout the paper, we  use the standard tensor notation, where for any  
$$A=\sum_{i,j,k,l=1}^N a_{ijkl} \ts e_{ij}\ot e_{kl}\,\in\, \ndo\CC^N\ot\ndo\CC^N$$
and  $r,s=1,\ldots , n$ with  $r\neq s$  and  $n\geqslant 2$      we denote by $A_{rs}$ the element of   $(\ndo\CC^N)^{\ot n}$, 
$$
A_{rs}=\sum_{i,j,k,l=1}^N a_{ijkl}\ts (e_{ij})_r (e_{kl})_s \quad\text{for}\quad
(e_{ij})_t = 1^{\ot (t-1)} \ot e_{ij} \ot 1^{\ot{(n-t)}}.
$$
For example, in \eqref{YBE} (resp. \eqref{unitrig}) we used this notation with   $n=3$ (resp. $n=2$).

Recall
the  formal Taylor Theorem, 
\beq\label{taylor}
a(z+z_0)=e^{z_0\frac{\partial}{\partial z} a(z)}=
\sum_{k=0}^\infty \frac{z_0^k}{k!} \frac{\partial^k}{\partial z^k} a(z)\quad\text{for}\quad a(z)\in V[[z^{\pm 1}]],
\eeq
where $V$ is a  vector space.
For example, by using \eqref{taylor}, we get the expansion
\beq\label{exp}
\frac{1}{1-e^{ah}x}=
\frac{1}{1-(x+(e^{ah}-1)x)}=
\sum_{k=0}^{\infty} \frac{(e^{ah}-1)^k x^k}{k!}\frac{\partial^k}{\partial x^k}\left(\frac{1}{1-x}\right), \quad a\in\CC,
\eeq
where $(1-x)^{-1}=\sum_{k\geqslant 0}x^k\in\CC[[x]]$.
Suppose $q$ is a formal parameter. By  \cite{FR}, there exists a unique power series 
$f_q (x)=1+\sum_{k=1}^{\infty} f_{q,k} x^k(1-x)^{-k} \in\CC(q)[[x]]$,
where  all $f_{q,k} (q-1)^{-k} \in\CC(q)$ are regular at $q=1$ for all $k\geqslant 1$ (see also \cite{KM}),
such that we have
$$
f_q (xq^{2N}) = f_q (x)\frac{\left(1-xq^2\right)\left(1-xq^{2N-2}\right)}{\left(1-x\right)\left(1-xq^{2N}\right)}.
$$
 Therefore, by applying the substitution $q=e^{h/2}$ to $f_q(x)$  
  we obtain the power series
\beq\label{f3}
f(x)=1+\sum_{k=1}^{\infty} f_{k} \left(\frac{x}{1-x}\right)^k\in\CC[[x,h]],\quad\text{where}\quad f_k \coloneqq (f_{q,k})\left|_{q=e^{h/2}}\right.\in h^k\CC[[h]] .
\eeq

We shall also need the {\em normalized $R$-matrix},
\beq\label{R}
R(x)= ( 1-e^{-h}x )^{-1}  f(x)\ts \R(x)\,\in\, \ndo\CC^N \ot \ndo\CC^N[[x,h]].
\eeq
It
possesses the {\em crossing symmetry properties} (cf. \cite{FR})
\beq\label{csym}
R(xe^{Nh})^{t_1} \ts D_1  \ts ( R(x)^{-1})^{t_1}=D_1\fand (R(x)^{-1})^{t_2}  \ts D_2  \ts R(xe^{Nh})^{t_2} = D_2,
\eeq
where  $t_i$ is the matrix transposition $e_{rs}\mapsto e_{sr}$ applied on the $i$-th tensor factor and  $D $ is the diagonal $N\times N$ matrix
$D=\diag\left(e^{ (N-1)h/2 },e^{ (N-3)h/2},\ldots ,e^{- (N-1)h/2} \right)$.
As  $(e^{ah}-1) x$ lies in $ xh\CC[[h]]$,    by     \eqref{exp} and \eqref{f3}   the normalization term in \eqref{R} equals
\beq\label{g}
g(x)\coloneqq   \frac{f(x)}{1-e^{-h}x}=\sum_{k=0}^{\infty} g_k \frac{x^k}{\left(1-x\right)^{k+1}},\quad\text{where}\quad 
g_k\in h^k\CC[[h]]\fand g_0 =1.
\eeq
Thus, the   $R$-matrix $R(x)$ can be also   regarded as an element of $\ndo\CC^N \ot \ndo\CC^N(x)[[h]]$.

Let $\CC_* (z_1,\ldots ,z_n)$ be
 the localization  of the ring of Taylor series
$\CC[[z_1,\ldots ,z_n]]$ at  $\CC[z_1,\ldots ,z_n]^{\times}$.
There exists a unique embedding $ \CC_* (z_1,\ldots ,z_n)\hookrightarrow\CC((z_1))\ldots ((z_n))$, which  naturally extends to   $\CC_* (z_1,\ldots ,z_n)[[h]]$, so that we obtain the map
$$
\iotan\colon \CC_* (z_1,\ldots ,z_n)[[h]]\to \CC((z_1))\ldots ((z_n))[[h]].
$$
As in \cite{KM}, by  applying the substitution $x=e^u$ and then the embedding $\iotau$ to the normalization term $g(x)$ from \eqref{g}, we get  $\iotau  g(e^u) \in \CC((u))[[h]]$.
  Finally, 
	by \cite[Prop. 1.2]{EK4} and \cite[Prop. 2.1]{KM},
	there exists a unique $\psi\in 1+h\CC[[h]]$ such that  the $R$-matrix
\beq\label{rplusg}
R(e^u)\coloneqq \psi \iotau g(e^u)   \R  (e^u)  \,\in\, \ndo\CC^N\ot\ndo\CC^N((u))[[h]]
\eeq
has the {\em unitarity property}
$ 
R_{12}(e^u) R_{21}(e^{-u}) =1
$ 
and 
the {\em crossing symmetry properties}
$$
R(e^{u+Nh})^{t_1}  D_1 ( R(e^u)^{-1})^{t_1}=D_1\fand (R(e^u)^{-1})^{t_2} D_2 R( e^{u+Nh})^{t_2} = D_2.
$$

Throughout the paper, whenever it is clear from the context, we      omit   the  embedding symbol $\iotaop$, e.g., we write   $g(e^u)$ instead of $\iotau g(e^u)$. In the multiple variable case, we   use the common expansion convention where the   embedding is determined by the order of the variables, e.g., if $\sigma$ is a permutation in the symmetric group $\mathfrak{S}_n$, then $g(e^{u_{\sigma_1} +\ldots +u_{\sigma_n}})$   denotes the series $\iotasigma  g(e^{u_{\sigma_1} +\ldots +u_{\sigma_n}})$ which belongs to $\CC((u_{\sigma_1}))\ldots ((u_{\sigma_n} ))[[h]]$.

\subsection{Quantum affine vertex algebra \texorpdfstring{$\Vcgl$}{Vc(glN)}}\label{subsection_12}

We follow \cite{EK3,EK4} to introduce the {\em quantized universal enveloping algebra} $\Ur$; see also \cite{FRT,RS}.
It is defined as an associative algebra over the ring $\CC[[h]]$ generated by  the elements $t_{ij}^{(-r)}$, where $i,j=1,\ldots ,N$ and $r=1,2,\ldots$, subject to the  defining relations given by
\beq\label{rtt}
R(e^{u-v})\ts T_{1}^+(u)\ts T_2^+ (v)=  T_2^+ (v)\ts T_{1}^+(u)\ts R(e^{u-v}).
\eeq
In \eqref{rtt},  we use subscripts to indicate the tensor factors in   
$\ndo\CC^N \ot \ndo\CC^N \ot\Ur$  and the matrix $T^+(u)  $ is given by
$$
T^+(u) =\sum_{i,j=1}^N e_{ij}\ot t_{ij}^+ (u),\qquad\text{where}\qquad t_{ij}^+(u)=\delta_{ij}-h\sum_{r=1}^{\infty}t_{ij}^{(-r)}u^{r-1}\in\Ur[[u]].
$$
It is worth noting that the algebra $\Ur$ is closely connected with the dual Yangian for $\mathfrak{gl}_N$; see \cite[Sect. 2]{KM}. Throughout the paper, we assume that $\Ur$ is $h$-adically completed. In addition,   the tensor products of  $h$-adically complete $\CC[[h]]$-modules  are understood as $h$-adically completed;
see
\cite[Chapter XVI]{Kas} for more information on the $h$-adic topology.

Let $P^h\in \ndo\CC^N \ot \ndo\CC^N [[h]] $ be the $h$-permutation operator,
$$
P^h = \sum_{i=1}^N e_{ii}\ot e_{ii} + e^{h/2}\sum_{\substack{i,j=1\\i> j}}^N e_{ij}\ot e_{ji} +e^{-h/2}\sum_{\substack{i,j=1\\i< j}}^N e_{ij}\ot e_{ji}.
$$
Consider the action of the symmetric group $\mathfrak{S}_N$ on the space $(\CC^N)^{\ot N}$  given by $\sigma_i\mapsto P_{\sigma_i}^h=P_{i\ts i+1}^h$ for $i=1,\ldots ,N-1$, where
$\sigma_i$ denotes the transposition $(i,i+1)$.  Let $A^{(N)}$ be the image of the normalized anti-symmetrizer under this action,  
$$
A^{(N)}=\frac{1}{N!}\sum_{\sigma\in\mathfrak{S}_N} \sgn\sigma \cdot P_\sigma^h .
$$
Consider the {\em quantum determinant} of the matrix $T^+(u)$,
\beq\label{defofqdet}
\qdet T^+ (u)=\tr_{1,\ldots ,N} \,A^{(N)}\ts T_1^+(u)\ldots T_N^+(u-(N-1)h)\ts D_1\ldots D_N \in \Ur [[u]],
\eeq
where the trace is taken over all $N$ copies of $\ndo\CC^N$. Its coefficients $\delta_r$,  given by
\beq\label{detkoef}
\qdet T^+ (u)=1-h\sum_{r\geqslant 0}\delta_r u^r,
\eeq
belong to the center of the algebra $\Ur$; see proof of \cite[Prop. 3.10]{KM}.

Our next goal is to recall the construction of the {\em quantum affine vertex algebra} $\Vcgl$ from \cite{EK5}.
First, for  
the families of variables $u = (u_1 ,...,u_n )$ and $v = (v_1 ,...,v_m )$ and $a\in\CC$, introduce the formal power series with coefficients in
$(\ndo\CC^N )^{\ot n} \ot (\ndo\CC^N )^{\ot m}$
 by
\begin{align}
&R_{nm}^{12}(e^{z+u-v+ah})= \prod_{i=1,\dots,n}^{\longrightarrow} 
\prod_{j=n+1,\ldots,n+m}^{\longleftarrow}   R_{ij} (e^{z+u_i-v_{j-n}+ah}),\label{rnm12exp}
\end{align}
where the arrows indicate the order of factors.
Also,   denote by $R_{nm}^{12}(e^{ u-v+ah})$ the analogous product of the $R$-matrices  $R_{ij} (e^{ u_i-v_{j-n}+ah})$.
Let $\vac$ be the unit in the algebra $\Ur$. From now on, we regard the coefficients of the matrix entries of   $T^+(u)$ as operators on $\Ur$.
Moreover, we indicate by  the parameter $h$ in the subscript    that the corresponding $\CC[[h]]$-module is $h$-adically completed, e.g., if  $V$ is a topologically free $\CC[[h]]$-module,   then  
$V((z ))_h$  is   the $\CC[[h]]$-module of all $a(z)=\sum_{r \in\ZZ}a_{ r}\ts z^{-r}\in V[[z^{\pm 1}]]$ such that $\lim_{r\to\infty} a_r =0$ in the $h$-adic topology, i.e., the $h$-adic completion of $V((z))$.
Let us recall \cite[Lemma 2.1]{EK5}.

\begin{lem}\label{lemma21}
For any $c\in\CC$ there exists a unique operator series
$ T^*(u) \in\ndo\CC^N \ot \om( \Ur ,\Ur ((u))_h )$ 
such that   for all $n\geqslant 0$ we have
\beq\label{tstar_formulaa}
T^{*}_{0} (u)\ts T_{1}^{+}(v_1)\ldots T_{n }^+(v_n)\vac  = R_{1n}^{12}(e^{u-v+hc/2})^{-1} T_{1}^{+}(v_1)\ldots T_{n }^+(v_n) R_{1n}^{12}(e^{u-v-hc/2})\vac. 
\eeq
\end{lem}

From now on, we denote the topologically free $\CC[[h]]$-module of $\Ur$ by $\Vcgl$, where the superscript $c$ determines the action of the series $T^*(u)$ from \eqref{tstar_formulaa}. We  shall use the following abbreviations for the   products of the operators $T^+(u)$ and $T^*(u)$:
\begin{align*}
T_{[n]}^{+}(u|z)=T_{1}^{+}(z+u_1)\ldots T_{n}^{+}(z+u_n)\fand T_{[n]}^*(u|z)=T_{1}^*(z+u_1)\ldots T_{n}^*(z+u_n),
\end{align*}
where $u=(u_1,\ldots ,u_n)$. Also, by omitting the variable $z$ we obtain the operator products  
$$T_{[n]}^{+}(u )=T_{1}^{+}(u_1)\ldots T_{n}^{+}(u_n)\fand T_{[n]}^*(u)=T_{1}^*(u_1)\ldots T_{n}^*(u_n) .$$
Let us  recall the Etingof--Kazhdan construction  \cite[Thm. 2.3]{EK5}; see also \cite[Sect. 1.4.1]{EK5} for the definition of quantum vertex algebra.
 
\begin{thm}\label{EK:qva}
For any $c\in \CC$
there exists a unique  quantum vertex algebra  structure on $\Vcgl$
  such that the vertex operator map $Y(\cdot, z)$ is given by
\beq\label{EK_Ymap}
Y\big(T_{[n]}^+ (u)\vac,z\big)=T_{[n]}^+ (u|z)\ts T_{[n]}^* (u|z+hc/2)^{-1}, 
\eeq
the vacuum vector is $\vac $
  and the braiding  map  is defined by the relation  
\begin{align}
&\mathcal{S}(z)\big(R_{nm}^{  12}(e^{z+u-v})^{-1} \ts  T_{[m]}^{+24}(v)  \ts  
R_{nm}^{  12}(e^{z+u-v-h  c})  \ts   T_{[n]}^{+13}(u)(\vac\otimes \vac) \big) \non \\
 =&\, T_{[n]}^{+13}(u)  \ts   R_{nm}^{  12}(e^{z+u-v+h  c})^{-1}  \ts  
 T_{[m]}^{+24}(v)  \ts   R_{nm}^{  12}(e^{z+u-v})(\vac\otimes \vac)  \label{braiding}
\end{align}
 on
$
(\ndo\mathbb{C}^{N})^{\otimes n} \otimes
(\ndo\mathbb{C}^{N})^{\otimes m}\otimes \Vcgl \ot \Vcgl $.
\end{thm}

Extend the notation   \eqref{rnm12exp} by
\begin{align}
&R_{nm}^{12}(ze^{u-v+ah})= \prod_{i=1,\dots,n}^{\longrightarrow} 
\prod_{j=n+1,\ldots,n+m}^{\longleftarrow}   R_{ij} (ze^{u_i-v_{j-n}+ah}),\label{rnm12expx}  
\end{align}
where  the $R$-matrices $R_{ij} (ze^{u_i-v_{j-n}+ah})$   are regarded as rational functions in   $z$.  
Finally, we recall \cite[Lemma 2.6]{K}, which gives us a multiplicative version of the braiding   \eqref{braiding}.

\begin{lem}\label{S_hat_lemma} 
There exists a unique $\CC[[h]]$-module map 
$\wht{\Sc}(z)\colon \Vcgl^{\ot 2}   \to  \Vcgl^{\ot 2}  \ot\CC(z)[[h]]$ 
such that   on
$
(\ndo\mathbb{C}^{N})^{\otimes n} \otimes
(\ndo\mathbb{C}^{N})^{\otimes m}\otimes \Vcgl \ot \Vcgl $ we  have
\begin{align}
&\wht{\Sc}(z)\big(R_{nm}^{  12}(ze^{ u-v})^{-1} \ts  T_{[m]}^{+24}(v)  \ts  
R_{nm}^{  12}(ze^{ u-v-h  c})  \ts   T_{[n]}^{+13}(u)(\vac\otimes \vac) \big) \non \\
 =&\, T_{[n]}^{+13}(u)  \ts   R_{nm}^{  12}(ze^{ u-v+h  c})^{-1}  \ts  
 T_{[m]}^{+24}(v)  \ts   R_{nm}^{  12}(ze^{ u-v})(\vac\otimes \vac) .\label{s-hat_def}
\end{align}
\end{lem}

\section{Multiplicative compatible pairs over \texorpdfstring{$\Vcgl$}{Vc(glN)}}\label{section_02}

In this section, we adapt the notion of compatible pair from \cite[Sect. 3]{BK}, so that it fits the setting of $\phi$-coordinated modules  \cite{Liphi}. Furthermore, we construct an example of such a pair over $\Vcgl$. Throughout the rest of the paper, the level $c\in\CC$ is arbitrary.

\begin{defn}\label{def_br_mult}
Let $V$ be a topologically free $\CC[[h]]$-module. The map
$$
\mu(z,x)\colon V\ot V\to V\ot V\ot\CC (z,x) [[h]]
$$
is said to be a {\em multiplicative  braiding map} if it satisfies the {\em quantum Yang--Baxter equation}
\begin{align}
&\mu_{12}(z_1,x z_2)\ts\mu_{13}(z_1z_2,x  )\ts\mu_{23}(z_2,x )=
\mu_{23}(z_2,x )\ts\mu_{13}(z_1z_2,x  )\ts\mu_{12}(z_1,x z_2), \label{ybe_muh}\\
\intertext{  the {\em unitarity condition}}
&\mu(1/z,x)\ts \mu_{21}(z,x/z)=\mu_{21}(z,x/z)\ts\mu(1/z,x) =1 \label{uni_muh}
\intertext{  and the {\em singularity constraint}}
&    \mu_{\pm}(z_1 /z_2,z_2)(u\ot v) \in V\ot V\ot\CC[z_1,z_1^{-1}]((z_2^{\pm 1}))[[h]]\text{ for all }u,v\in V,\label{sing_muh}\\
\intertext{where}
&   \mu_{\pm}(z_1 /z_2,z_2)(u\ot v) 
=\left(\iotazepm  \mu(z ,x)(u\ot v)\right)\big|_{z=z_1 /z_2,\ts x=z_2}.\non
\end{align}
\end{defn}

\begin{rem}\label{remark_32}
In contrast with \cite[Def. 3.1]{BK}, where the braiding map is required to satisfy the additive versions of the   Yang--Baxter equation and the unitarity condition,
the   map in Definition \ref{def_br_mult} possesses multiplicative counterparts of these properties. To emphasize this distinction, we refer to it as a {\em multiplicative} braiding map in the definition.
However, aside from  the Etingof--Kazhdan braiding \eqref{braiding},
throughout this paper we consider   multiplicative braiding maps only, so we suppress the term multiplicative from now on.
\end{rem}

\begin{rem}
Clearly, the requirements \eqref{ybe_muh} and \eqref{uni_muh} imposed on the braiding map are motivated by the corresponding properties from the definition of quantum vertex algebra; see \cite[Sect. 1.4.1]{EK5}. As for the singularity constraint \eqref{sing_muh}, we explain its role in Remarks  \ref{rem_sing_1} and  \ref{rem_sing_2} below. At this point, we only remark that the condition \eqref{sing_muh} is nontrivial as, e.g., it does not hold for the function
$
m(z,x)= (1-zx)^{-1} \in\CC(z,x)\subset\CC(z,x)[[h]]
$.
\end{rem}

\begin{ex}\label{ex_34}
The map $\wht{\Sc} =\wht{\Sc} (z)$, as  defined by \eqref{s-hat_def}, 
is a braiding map, which is constant with respect to the second variable $x$. Indeed, this is due to the fact that the  additive counterparts of the properties \eqref{ybe_muh} and \eqref{uni_muh} hold for the Etingof--Kazhdan braiding $\Sc$ defined by \eqref{braiding}; see \cite{EK5} for more details. 
\end{ex}
 
To construct another example of braiding map, we shall need the notation
\begin{align}
&\RR_{nm}^{12\ts t}(1/x^2 z e^{ u+v})= \prod_{i=1,\dots,n}^{\longrightarrow}
\prod_{j=n+1,\ldots,n+m}^{\longrightarrow}  \R_{ij}(1/x^2 z e^{u_i +v_{j-n}})^{t_i},\label{r1rrnm}\\
&\RR_{nm}^{12\ts *}(1/x^2 z e^{ u+v}) = \prod_{i=1,\dots,n}^{\longleftarrow}
\prod_{j=n+1,\ldots,n+m}^{\longleftarrow}  \left(\R_{ij}(1/x^2 z e^{u_i +v_{j-n}})^{-1}\right)^{t_i}, \label{r2rrnm}
\end{align}
 where $u=(u_1,\ldots ,u_n)$, $v=(v_1,\ldots ,v_m)$,  $z$ and $x$ are single variables and, as before, the arrows   indicate the order of factors. For example, we have
$$
\RR_{23}^{12\ts t}(1/x^2 z e^{ u+v})=\R_{13}^t \R_{14}^t \R_{15}^t \R_{23}^t \R_{24}^t \R_{25}^t,\quad 
\RR_{23}^{12\ts *}(1/x^2 z e^{ u+v})=\R_{25}^*\R_{24}^*\R_{23}^*\R_{15}^*\R_{14}^* \R_{13}^*,
$$
for $\R_{ij}^t =\R_{ij}(1/x^2 z e^{u_i +v_{j-n}})^{t_i}$, $\R_{ij}^* = \left(\R_{ij}(1/x^2 z e^{u_i +v_{j-n}})^{-1}\right)^{t_i}$.
 In addition, we shall write
$\R_{nm}^{12}(z e^{ u-v+ah})$ for the   product  defined as in \eqref{rnm12expx}, with the $R$-matrices $R_{ij} (ze^{u_i-v_{j-n}+ah})$ now replaced by $\R_{ij} (ze^{u_i-v_{j-n}+ah})$. Finally,  denote by 
``$\cdotlr$'' (resp. ``$\cdotrl$'')   the standard multiplication in $\ndo\CC^N\ot(\ndo\CC^N)^{\text{op}}$ (resp. $(\ndo\CC^N)^{\text{op}}\ot \ndo\CC^N $), where $A^{\text{op}}$  is the opposite algebra of $A$. One easily checks that the $R$-matrix products in \eqref{r1rrnm} and \eqref{r2rrnm} exhibit the property, which we often use in the rest of the paper,
\beq\label{csym_cons}
\RR_{nm}^{12\ts t}(1/x^2 z e^{ u+v})\cdotlr\RR_{nm}^{12\ts *}(1/x^2 z e^{ u+v}) 
=\RR_{nm}^{12\ts t}(1/x^2 z e^{ u+v})\cdotrl\RR_{nm}^{12\ts *}(1/x^2 z e^{ u+v}) 
=1,
\eeq
where ``$\cdotlr$'' (resp. ``$\cdotrl$'')  is the standard multiplication in $(\ndo\CC^N)^{\ot n}\ot((\ndo\CC^N)^{\ot m})^{\text{op}}$ (resp. $((\ndo\CC^N)^{\ot n})^{\text{op}}\ot (\ndo\CC^N)^{\ot m} $).

\begin{pro}\label{pro_muhat}
There exists a unique    braiding map $\mu$ on $\Vcgl$
such that  we have
\begin{align}
\mu (z,x)\left(T_{[n]}^{+13}(u)T_{[m]}^{+24}(v)(\vac\ot\vac)\right)
=
\RR_{nm}^{12\ts *}(1/x^2 z e^{ u+v})
\cdotlr\left(\R_{nm}^{12}(ze^{u-v})\ts T_{[n]}^{+13}(u)\ts \right.&\non\\
\left. \times\RR_{nm}^{12\ts t}(1/x^2 z e^{ u+v}) \ts T_{[m]}^{+24}(v)\ts \R_{nm}^{12}(ze^{u-v})^{-1}(\vac\ot\vac)\right)&.\label{muh_formula}
\end{align}
\end{pro}

\begin{prf}
The fact that \eqref{muh_formula} defines a $\CC[[h]]$-module map follows from the  Yang--Baxter equation \eqref{YBE} for the $R$-matrix $\R(x)$. To prove it, it   suffices to show that the assignments 
\begin{align}
&\mu^{(1)} (z,x)\left(T_{[n]}^{+13}(u)T_{[m]}^{+24}(v)(\vac\ot\vac)\right)
=
 T_{[n]}^{+13}(u)\ts  \RR_{nm}^{12\ts t}(1/x^2 z e^{ u+v}) \ts T_{[m]}^{+24}(v) (\vac\ot\vac), \label{tmp_muu_1}
\\
&\mu^{(2)} (z,x)\left(T_{[n]}^{+13}(u)T_{[m]}^{+24}(v)(\vac\ot\vac)\right)
=
 T_{[n]}^{+13}(u)\ts T_{[m]}^{+24}(v)\ts\R_{nm}^{12}(ze^{u-v})^{-1} (\vac\ot\vac) ,\label{tmp_muu_2}
\\
&\mu^{(3)} (z,x)\left(T_{[n]}^{+13}(u)T_{[m]}^{+24}(v)(\vac\ot\vac)\right)
=
 \R_{nm}^{12}(ze^{u-v})\ts T_{[n]}^{+13}(u)\ts T_{[m]}^{+24}(v)(\vac\ot\vac), \label{tmp_muu_3}
\\
&\mu^{(4)} (z,x)\left(T_{[n]}^{+13}(u)T_{[m]}^{+24}(v)(\vac\ot\vac)\right)
=
T_{[m]}^{+24}(v)\ts
\RR_{nm}^{12\ts *}(1/x^2 z e^{ u+v})\ts
   T_{[n]}^{+13}(u)(\vac\ot\vac) \label{tmp_muu_4}
\end{align}
define $\CC[[h]]$-module maps $\mu^{(i)}$, $i=1,2,3,4$, as   $\mu$
is then equal to their composition  $ \mu^{(1)}\mu^{(2)}\mu^{(3)}\mu^{(4)}$.  Consider the expressions with coefficients in
$$
(\ndo\CC^N)^{\ot n} \ot (\ndo\CC^N)^{\ot m}\ot \Vcgl\ot\Vcgl,
$$
which are given by
\begin{align*}
&\left(R_{i\ts i+1}(e^{u_i -u_{i+1}})\ts T_{[n]}^{+13} (u)-
P_{i\ts i+1}\ts T_{[n]}^{+13} (u)^\prime \ts P_{i\ts i+1}\ts R_{i\ts i+1}(e^{u_i -u_{i+1}})\right)  T_{[m]}^{+24} (v)(\vac\ot\vac),  \\
&T_{[n]}^{+13} (u) \left(R_{j\ts j+1}(e^{v_{j-n} -v_{j-n+1}})\ts T_{[m]}^{+24} (v)-
P_{j\ts j+1}\ts T_{[m]}^{+24} (v)^\prime \ts P_{j\ts j+1}\ts R_{j\ts j+1}(e^{v_{j-n} -v_{j-n+1}})\right)(\vac\ot\vac) , 
\end{align*}
where $i=1,\ldots ,n-1$, $j=n+1,\ldots ,n+m-1$, the prime in $T_{[n]}^{+13} (u)^\prime$ (resp. $T_{[m]}^{+24} (v)^\prime$) indicates that the variables $u_i$ and $u_{i+1}$ (resp. $v_{j-n}$ and $v_{j-n+1}$) are swapped and $P_{a\ts a+1}$ stands for the action of the permutation operator
\beq\label{permutation_operator}
P=\sum_{i,j=1}^N e_{ij}\ot e_{ji}\in \ndo\CC^N \ot \ndo\CC^N
\eeq
on the tensor factors $a$ and $a+1$. One checks  directly that for all $n$ and $m$ the given expressions belong to the kernels of \eqref{tmp_muu_1}--\eqref{tmp_muu_4}, thus showing that $\mu$ is well-defined. Furthermore, it is clear that \eqref{muh_formula} uniquely determines the map $\mu$, as the monomials in   generators $t_{ij}^{(-r)}$, $i,j=1,\ldots , N$ and $r=1,2,\ldots ,$ form an $h$-adically dense subset of $\Ur$.

The Yang--Baxter equation \eqref{ybe_muh} for the map $\mu$ is again proved by a direct calculation which relies on \eqref{YBE}. We omit its proof as it goes in parallel with \cite[Lemma 4.6]{BK}. The proof of unitarity condition \eqref{uni_muh} is also straightforward. It relies on the   expression
\begin{align*}
\mu_{43} (z,x)\left(T_{[n]}^{+13}(u)T_{[m]}^{+24}(v)(\vac\ot\vac)\right)
=
\RR_{nm}^{12\ts *}(1/x^2 z e^{ u+v})
\cdotrl\left(\R_{nm}^{21}(ze^{v-u})\ts T_{[m]}^{+24}(v)\ts \right.&\non\\
\left. \times\RR_{nm}^{12\ts t}(1/x^2 z e^{ u+v}) \ts T_{[n]}^{+13}(u)\ts \R_{nm}^{21}(ze^{v-u})^{-1}(\vac\ot\vac)\right)&,
\end{align*}
 derived from \eqref{muh_formula}, where
$$
\R_{nm}^{21}(ze^{v-u})= \prod_{i=1,\dots,n}^{\longleftarrow} 
\prod_{j=n+1,\ldots,n+m}^{\longrightarrow}   \R_{ji} (ze^{v_{j-n}-u_i}),
$$
and, in addition, on the properties \eqref{unitrig} and
$\R(x)^{t_1}=\R_{21}(x)^{t_2}$
of the $R$-matrix \eqref{Rbar}.
Finally, the singularity constraint \eqref{sing_muh} is evident from the form of the $R$-matrices which appear in the defining expression \eqref{muh_formula} for the map $\mu$.
\end{prf}

\begin{defn}\label{def_int_mult}
Let $V$ be a topologically free $\CC[[h]]$-module.
A    braiding map $\mu \colon V\ot V\to V\ot V\ot\CC    (z,x)[[ h]]$ is said to be {\em compatible} if there  exists a  map  
$$
\nu  (z,x)\colon V\ot V\to V\ot V\ot\CC     (z,x) [[h]]
$$
which 
possesses the following properties. 
 The map $\nu $
is invertible in the sense that there exists a map
$  
\nu^{-1} (z,x)\colon V\ot V\to V\ot V\ot\CC    (z,x) [[h]]
$  
such that we have
\beq\label{def_int_mult_prop1}
\nu  (z,x)\ts\nu^{-1} (z,x)=\nu^{-1} (z,x)\ts \nu  (z,x)=1,
\eeq
and the map 
$ 
 \sigma  (z,x)\colon V\ot V\to V\ot V\ot\CC   (z,x)[[ h]]
$ 
defined by
$$
 \sigma (z,x)= \nu (z,x)\ts \mu (z,x)\ts \nu_{21}^{-1}(1/z,zx)
$$
is again a   braiding map (as defined by Definition \ref{def_br_mult}).  
The pair $(\mu ,\nu )$ is then said to be a {\em multiplicative compatible pair}.
\end{defn}

\begin{rem}
We shall usually omit the term {\em multiplicative} and refer to a pair $(\mu,\nu)$ satisfying Definition \ref{def_int_mult} more briefly as a compatible pair; recall Remark \ref{remark_32}.
\end{rem}

\begin{ex}
Clearly, for any braiding map $\mu \colon V\ot V\to V\ot V\ot\CC    (z,x)[[ h]]$, the pair $(\mu,\iota)$, where $\iota$ is the inclusion $V\ot V\hookrightarrow V\ot V\ot\CC(z,x)[[h]]$, is a   compatible pair. In particular, by Example \ref{ex_34}, the pair $(\wht{\Sc},\iota)$ is compatible.
\end{ex}

In the following proposition, we give another example of a compatible pair.

\begin{pro}\label{nuhatprop}
There exists a unique map
$$
\nu (z,x)\colon\Vcgl\ot\Vcgl\to \Vcgl\ot\Vcgl\ot\CC( z,x )[[ h]]
$$
such that we have
\begin{align}
&\nu (z,x)\left(T_{[n]}^{+13}(u)T_{[m]}^{+24}(v)(\vac\ot\vac)\right)\non\\
=&\,
\RR_{nm}^{12\ts *}(1/x^2 z  e^{ u+v})
\cdotrl
\left(
 T_{[n]}^{+13}(u)\ts
R_{nm}^{12}(z e^{ u-v+hc})^{-1}\ts
T_{[m]}^{+24}(v)\ts R_{nm}^{12}(z e^{ u-v})(\vac\ot\vac)\right).\label{nuhatformula}
\end{align}
The pair $(\mu ,\nu )$ is compatible and it satisfies
\beq\label{s_hat_id}
\wht{\Sc}(z)=\nu (z,x)\ts \mu (z,x)\ts \nu^{-1}_{21}(1/z,z x).
\eeq
\end{pro}

\begin{prfn}
One can check that $\nu$ is well-defined  by arguing as in the corresponding part of the proof of Proposition \ref{pro_muhat}.
Next, by using the second crossing symmetry property in \eqref{csym}, one proves that the map 
$\nu^{-1} (z,x)\colon \Vcgl\ot  \Vcgl\to  \Vcgl\ot  \Vcgl\ot\CC    (z,x) [[h]]$
defined by
\begin{align*}
 \nu^{-1} (z,x)\left(T_{[n]}^{+13}(u)T_{[m]}^{+24}(v)(\vac\ot\vac)\right)=\left((D_{[m]}^2)^{-1}\ts R_{nm}^{12}(z e^{ u-v+h(c+N)}) \ts D_{[m]}^2\right)&\\
\cdotrl
\left(
 T_{[n]}^{+13}(u)\ts
\RR_{nm}^{12\ts t}(1/x^2 z  e^{ u+v})\ts
T_{[m]}^{+24}(v)\ts R_{nm}^{12}(z e^{ u-v})^{-1}(\vac\ot\vac)\right)&, 
\end{align*}
where $D_{[m]}^2=D_{n+1}\ldots D_{n+m}$,
satisfies the identities in \eqref{def_int_mult_prop1}.
Finally, by Example \ref{ex_34}, we already know that $\wht{\Sc}$ is a braiding map. Hence, to finish the proof, it is sufficient to verify the identity \eqref{s_hat_id}. This is done by employing the explicit expressions for the   maps $\wht{\Sc}$, $\mu$ and $\nu$, as given by \eqref{s-hat_def}, \eqref{muh_formula} and \eqref{nuhatformula}, along with the formula
\begin{align*}
&\nu_{43} (1/z,zx)\left(T_{[n]}^{+13}(u)T_{[m]}^{+24}(v)(\vac\ot\vac)\right) \\
=&\,
\RR_{nm}^{12\ts *}(1/x^2 z  e^{ u+v})
\cdotlr
\left(
T_{[m]}^{+24}(v)\ts
R_{nm}^{12}(z e^{ u-v-hc}) \ts
 T_{[n]}^{+13}(u)\ts
 R_{nm}^{12}(z e^{ u-v})^{-1}(\vac\ot\vac)\right), 
\end{align*}
which is derived from \eqref{nuhatformula}.
More specifically,  one   shows that, equivalently, the images of 
$T_{[n]}^{+13}(u)T_{[m]}^{+24}(v)(\vac\ot\vac)$ 
under the maps 
$\wht{\Sc}(z)\nu_{21}(1/z,z x)$ and $\nu (z,x) \mu (z,x)$ coincide. Indeed, both are equal to
$$\pushQED{\qed} 
\RR_{nm}^{12\ts *}(1/x^2 z  e^{ u+v}) 
\cdotlr
\left(
 R_{nm}^{  12}(ze^{ u-v}) \ts T_{[n]}^{+13}(u)  \ts   R_{nm}^{  12}(ze^{ u-v+h  c})^{-1}  \ts  
 T_{[m]}^{+24}(v)   (\vac\otimes \vac)
\right).\eqno\qedhere\popQED
$$
\end{prfn}

Let $(\mu,\nu)$ be a compatible pair over $V$. We associate to $\nu=\nu(z,x)$ the map
$$
\nu(e^z,x) \colon V\ot V\to V\ot V\ot \CC_*(z,x)[[h]]
$$
defined by $\nu(e^z,x)=\nu(y,x)\big|_{y=e^z}$.
For example, the map $\nu $ given by \eqref{nuhatformula} yields
\begin{align}
&\nu (e^z,x)\left(T_{[n]}^{+13}(u)T_{[m]}^{+24}(v)(\vac\ot\vac)\right)\non\\
=&\,
\RR_{nm}^{12\ts *}(1/x^2   e^{ z+u+v}) 
\cdotrl
\left(
 T_{[n]}^{+13}(u)\ts
R_{nm}^{12}(  e^{ z+u-v+hc})^{-1}\ts
T_{[m]}^{+24}(v)\ts R_{nm}^{12}(  e^{z+ u-v})(\vac\ot\vac)\right),\label{nuexpformula}
\end{align}
where  
$$
\RR_{nm}^{12\ts *}(1/x^2   e^{ z+u+v}) = \prod_{i=1,\dots,n}^{\longleftarrow}
\prod_{j=n+1,\ldots,n+m}^{\longleftarrow}  \left(R_{ij}(1/x^2   e^{z+u_i +v_{j-n}})^{-1}\right)^{t_i}. 
$$
Suppose, in addition, that $V$ is a quantum vertex algebra with the vertex operator map
$Y(\cdot,z)$. In the following sections, we shall  need
 the maps
$$
Y^{\nu }_{\pm}(\cdot,z,x)\colon V\ot V\to V((z))_h \ot \CC((x^{\pm 1}))((z))[[h]]
$$
 defined by
$$
Y^{\nu }_{\pm}(u,z,x)v= Y(z )\ts  \iotazepmrev   \nu  ( e^{z },x)(u\ot v)\quad\text{for all }u,v\in V.
$$
For example, by using \eqref{nuexpformula}, we find that, for   the vertex operator map \eqref{EK_Ymap}, we have
\beq\label{ynu_expl}
Y_{\pm}^{\nu} ( T_{[n]}^{+13}(u)\vac,z,x)  T_{[m]}^{+24}(v)\vac
=
\iotazepmrev \RR_{nm}^{12\ts *}(1/x^2   e^{ z+u+v}) 
\cdotrl
\left(
 T_{[n]}^{+13}(u|z)\ts
T_{[m]}^{+23}(v)\vac \right).
\eeq

We finish this section by   showing that the coefficients of the quantum determinant   \eqref{detkoef} give rise to the fixed points
of the
braiding map  \eqref{muh_formula}.
\begin{pro}\label{mufixes}
For any $w\in\Vcgl$ and $r=1,2,\ldots$ we have 
$$\mu (z,x) (w\ot \delta_r) =w\ot \delta_r.$$
\end{pro}

\begin{prf}
A particular case of the fusion procedure for the   $R$-matrix $\R(z)$, going back to \cite{C}, implies that 
for $u=(u_1,\ldots ,u_n)$ and $v=(v_0,v_0 -h ,\ldots ,v_0  -(N-1)h )$
there exist
$\alpha^{\pm 1},\beta^{\pm 1}$ in $ \CC(z,x)[[v_0,u_1,\ldots ,u_n,h]]$ 
such that  we have
\begin{alignat}{2}
&A^{(N)}\ts\RR_{nN}^{12\ts t}(1/x^2 z e^{ u+v})=\alpha  A^{(N)},\quad
&&A^{(N)}\ts\RR_{nN}^{12\ts *}(1/x^2 z e^{ u+v}) =\alpha^{-1}  A^{(N)},\label{alignat1}\\
&A^{(N)}\ts \R_{nN}^{12}(ze^{u-v})^{\pm 1} =\beta^{\pm 1}  A^{(N)},
&&A^{(N)}\ts   T_{[N]}^+(v )  =\cev{T}_{[N]}^+(v ) \ts A^{(N)}.\label{alignat2}
\end{alignat}
In fact, one can compute the explicit expressions for $\alpha^{\pm 1}$ and $\beta^{\pm 1}$, but we do not need them in this proof.
In \eqref{alignat1} and in the first equality in \eqref{alignat2}, the anti-symmetrizer  $A^{(N)}$ is applied on the tensor factors $n+1,\ldots ,n+N$ of $(\ndo\CC^N)^{\ot (n+N)}$. In the second equality  in \eqref{alignat2}, we use the left arrow on the top of the symbol to indicate that its factors are given in the opposite order, i.e.,  $\cev{T}_{[N]}^+(v )= T_N^+(v_0  -(N-1)h )\ldots T_{1}^+(v_0)$. The proposition can be now verified  by   usual arguments which rely  on \eqref{alignat1} and \eqref{alignat2}; see, e.g., the proof of \cite[Lemma 4.13]{BK}. However, we provide the proof details for completeness. Let us  apply  the braiding map $\mu $, as given by \eqref{muh_formula}, on the expression $T_{[n]}^{+13}(u)  A^{(N)}T_{[N]}^{+24}(v)(\vac\ot\vac)D_{[N]}^{2}$, where the anti-symmetrizer is again applied on the tensor factors $n+1,\ldots ,n+N$ and $D_{[N]}^{2}=D_{n+1}\ldots D_{n+N}$. Clearly, given the definition of the quantum determinant \eqref{defofqdet}, it suffices to show that the aforementioned expression is fixed by  $\mu$.
To simplify the notation, we introduce some  abbreviations for the $R$-matrix products,
$$
\RR_{nN}^{12\ts   *} =\RR_{nN}^{12\ts *}(1/x^2 z e^{ u+v}) ,\quad
\RR_{nN}^{12\ts t} =\RR_{nN}^{12\ts t}(1/x^2 z e^{ u+v}),\quad
(\R_{nm}^{12\ts })^{\pm 1}=\R_{nm}^{12}(ze^{u-v})^{\pm 1}.
$$
By employing  
\eqref{muh_formula}
and the identities given by \eqref{alignat1} and \eqref{alignat2}, 
we find that the image of $T_{[n]}^{+13}(u)  A^{(N)}T_{[N]}^{+24}(v)(\vac\ot\vac)D_{[N]}^{2}$ under the braiding map $\mu=\mu(z,x)$ equals
\begin{align*}
 &\, A^{(N)}\ts \RR_{nm}^{12\ts   *} 
\cdotlr\left(\R_{nm}^{12 } \ts T_{[n]}^{+13}(u)\ts \RR_{nm}^{12\ts t}  \ts T_{[N]}^{+24}(v)\ts (\R_{nm}^{12 })^{-1}(\vac\ot\vac) \right) D_{[N]}^{2}\\
=& \, \RR_{nm}^{12\ts  *} 
\cdotlr\left(A^{(N)}\ts\R_{nm}^{12 } \ts T_{[n]}^{+13}(u)\ts \RR_{nm}^{12\ts t}  \ts T_{[N]}^{+24}(v)\ts (\R_{nm}^{12 })^{-1} (\vac\ot\vac)\right)D_{[N]}^{2}\\
=&\, \beta \ts \RR_{nm}^{12\ts   *} 
\cdotlr\left( T_{[n]}^{+13}(u) A^{(N)}\ts\RR_{nm}^{12\ts t}  \ts T_{[N]}^{+24}(v)\ts (\R_{nm}^{12 })^{-1} (\vac\ot\vac)\right)D_{[N]}^{2}\\
=&\, \alpha\ts\beta \ts \RR_{nm}^{12\ts   *} 
\cdotlr\left( T_{[n]}^{+13}(u) A^{(N)}\ts  T_{[N]}^{+24}(v)\ts (\R_{nm}^{12 })^{-1} (\vac\ot\vac)\right)D_{[N]}^{2}\\
=&\, \alpha\ts\beta \ts \RR_{nm}^{12\ts   *} 
\cdotlr\left( T_{[n]}^{+13}(u)  \cev{T}_{[N]}^{+24}(v)\ts A^{(N)}\ts (\R_{nm}^{12 })^{-1} (\vac\ot\vac)\right)D_{[N]}^{2}\\
=&\, \alpha\ts  \RR_{nm}^{12\ts   *} 
\cdotlr\left( T_{[n]}^{+13}(u)  \cev{T}_{[N]}^{+24}(v)\ts A^{(N)}(\vac\ot\vac) \right)D_{[N]}^{2}\\
=& \,  T_{[n]}^{+13}(u)\ts  \cev{T}_{[N]}^{+24}(v)\ts A^{(N)} (\vac\ot\vac) D_{[N]}^{2} 
= T_{[n]}^{+13}(u)\ts   A^{(N)}\ts T_{[N]}^{+24}(v)(\vac\ot\vac)D_{[N]}^{2},
\end{align*}
as required.
\end{prf}

\section{Twisted \texorpdfstring{$h$}{h}-Yangian as a truncated deformed \texorpdfstring{$\phi$}{phi}-coordinated \texorpdfstring{$\Vcgl$}{Vc(glN)}-module}\label{section_03}

In this section, we introduce a $\phi$-coordinated analogue of the notion of deformed module   \cite[Def. 3.5]{BK}. In general,     {\em $\phi$-coordinated modules} can be considered for any   {\em associate $\phi=\phi(z_0,z_2)\in\CC((z_0))[[z_2]]$ of a one-dimensional additive formal group}; see \cite{Liphi} for more details.  The focus of this paper is   on   associating  quantum vertex algebra theory to a particular class of quantum algebras,  twisted $h$-Yangians. Thus, in comparison with \cite{BK}, we slightly specialize the next definition, so that it better fits this goal. In particular, we present it in the case of   associate $\phi(z_2,z_0)=z_2 e^{z_0}$, which is used throughout the paper.

\begin{defn}\label{defn_deformed_top}
Let $(V,Y,\vac,\Sc)$  be a    quantum vertex algebra,  $(\mu ,\nu )$     a
   compatible pair on $V$ and $W$   a topologically free $\CC[[h]]$-module equipped with  a $\CC[[h]]$-module map
\begin{align}
Y_W(\cdot ,z) \colon V\ot W&\to W[[z^{-1}]],\label{labe_defn_deformed_top} \\
v\ot w&\mapsto Y_W(z)(v\ot w)=Y_W(v,z)w=\sum_{r\in\mathbb{Z}_{\geqslant 0}} v_{r-1} w \ts z^{-r }.\non
\end{align}
A pair $(W,Y_W)$ is said to be a {\em  truncated\footnote{The term {\em truncated} indicates that the module map $Y_W(\cdot ,z)$ is  truncated from above; see \eqref{labe_defn_deformed_top}.}  $(\mu ,\nu )$-deformed $\phi$-coordinated $V$-module} if the   map $Y_W(\cdot, z)$ satisfies the {\em vacuum property}
\begin{align}
&Y_W(\vac, z)w=w\quad\text{for all }w\in W,\non\\
\intertext{the   {\em $\nu $-associativity}:
for any   $u,v\in V$   
we have}
 &   Y_W(u,z_2 e^{z_0})\ts Y_W(v,z_2)  
=    Y_W\big(Y^{\nu }_{-}(u,z_0,z_2)v,z_2\big)   , \non
\intertext{and the {\em $\mu $-commutativity}: for any $u,v\in V$ and $w\in W$   we have  }
&   Y_W(z_1)\big(1\otimes Y_W(z_2)\big)\big(  \mu_{-}(  z_1/z_2,z_2)(u\otimes v)\otimes w\big)= Y_W(v,z_2)  Y_W(u,z_1)  w    \label{ficommprop} 
.
\end{align}
\end{defn}

\begin{rem}\label{rem_sing_1}
Regarding Definition \ref{defn_deformed_top}, we remark that the singularity constraint \eqref{sing_muh}  ensures that the left-hand side of the $\mu$-commutativity property  \eqref{ficommprop} is well-defined. 
\end{rem}

Our next goal is to construct an example of truncated $(\mu ,\nu )$-deformed $\phi$-coordinated $\Vcgl$-module by using the structure of orthogonal twisted $h$-Yangian. First, we follow Molev, Ragoucy and Sorba \cite[Sect. 3.1]{MRS} to recall this algebra. In contrast with its original definition, which is given over the complex field, we introduce it  over the ring $\CC[[h]]$, so that it is in tune with the quantum vertex algebraic setting. The {\em twisted $h$-Yangian} $\Yo$ is the $h$-adically completed algebra generated by the elements $s_{ij}^{(r)}$, where $i,j=1,\ldots ,N$ and $r=0,1,\ldots ,$  subject to  the defining relations 
\begin{align}
&\R(x/y)\ts S_1(x) \ts  \R(1/xy)^{t_1}\ts S_2(y)
=S_2(y)\ts  \R(1/xy)^{t_1}\ts S_1(x) \ts \R(x/y),\label{refl_mult}\\
&s_{ij}^{(0)}=0\text{ for }1\leqslant i<j\leqslant N \fand
s_{ii}^{(0)}=1\text{ for }i=1,\ldots ,  N.\label{refl_mult_2nd}
\end{align}
The matrix $S(u)$ is defined by
$$
S(u)=\sum_{i,j=1}^N e_{ij}\ot s_{ij}(u),
$$ 
where for all $i,j=1,\ldots ,N$ its matrix entries $s_{ij}(u)$ are given by
$$s_{ii}(u)=1+h\sum_{r\geqslant 1} s_{ii}^{(r)}u^{-r} 
\Fand
s_{ij}(u)=h\sum_{r\geqslant 0} s_{ij}^{(r)}u^{-r}\quad \text{if}\quad i\neq j.
$$

The  {\em reflection equation} \eqref{refl_mult}     can be  generalized as follows. For the families of variables  $x=(x_1,\ldots ,x_n)$ and $y=(y_1,\ldots ,y_m)$, introduce the notation
\begin{align}
& \R_{nm}^{12 }(x/y) = \prod_{i=1,\dots,n}^{\longrightarrow}
\prod_{j=n+1,\ldots,n+m}^{\longleftarrow}   \R_{ij}(x_i / y_j) , \non\\
&\RR_{nm}^{12\ts t}(1/xy) = \prod_{i=1,\dots,n}^{\longrightarrow}
\prod_{j=n+1,\ldots,n+m}^{\longrightarrow}   \R_{ij}(1/x_i y_j) ^{t_i},\non\\
&\RR_{nm}^{12\ts *}(1/xy)  = \prod_{i=1,\dots,n}^{\longleftarrow}
\prod_{j=n+1,\ldots,n+m}^{\longleftarrow}   \left(\R_{ij}(1/x_i y_j)^{-1}\right) ^{t_i},\non\\
&S_{[n]}(x)=
  \prod_{i=1,\dots,n}^{\longrightarrow}
 S_i(x_i)\ts \R_{i\ts i+1}(1/x_i x_{i+1})^{t_i}\ldots  \R_{i  n}(1/x_i x_{n})^{t_i}.\label{snixovi} 
\end{align}
By combining  the quantum  Yang--Baxter equation \eqref{YBE} and the reflection equation \eqref{refl_mult}, one obtains the  general identity
\beq\label{rmatrixprod3refl}
\R_{nm}^{12}(x/y)\ts S^{13}_{[n]}(x) \ts  \RR_{nm}^{12\ts t}(1/xy) \ts S_{[m]}^{23}(y)
=\ts S_{[m]}^{23}(y)\ts  \RR_{nm}^{12\ts t}(1/xy) \ts S^{13}_{[n]}(x) \ts  \R_{nm}^{12}(x/y).
\eeq
Furthermore, it is worth noting the following decomposition of \eqref{snixovi}:
\beq\label{rmatrixprod4}
S_{[n+m]}(x,y)=S^{13}_{[n]}(x) \ts  \RR_{nm}^{12\ts t}(1/xy) \ts S_{[m]}^{23}(y).
\eeq

From now on, let $(\mu , \nu )$ be the compatible pair established by Proposition \ref{nuhatprop}.

\begin{thm}\label{thm_hyang} 
There exists a unique structure of truncated $(\mu , \nu )$-deformed $\phi$-coordinated $\Vcgl$-module on   $\Yo$   such that the module map is given by
\beq\label{yyoformula}
Y_{\Yo}(T_{[n]}^+(u_1,\ldots ,u_n)\vac,z) =  S_{[n]}(ze^{u_1},\ldots , ze^{u_n}).
\eeq
\end{thm}

\begin{prf}
We omit the proof as the theorem   is verified by the   arguments analogous to those from the proof of Theorem \ref{thm_hyang_ext} below. 
In fact, 
 the proofs of $\nu$-associativity and $\mu$-commutativity  from Definition \ref{defn_deformed_top} are   simpler than the proofs of their counterparts from Definition \ref{def_def_mod_phi}, which are  presented below. Indeed, due to the truncation condition \eqref{labe_defn_deformed_top}   imposed by Definition \ref{defn_deformed_top}, to adapt the proof of  Theorem \ref{thm_hyang_ext} to this setting, one can   choose for the polynomials $p$ and $q$ from its proof to be constant polynomial $1$.
\end{prf}

The next corollary is another  consequence of the proof of Theorem \ref{thm_hyang_ext}.

\begin{kor}\label{coro_conv_yo}
Let $W$ be a topologically free $\CC[[h]]$-module and also a $\Yo$-module. Denote the action of $S(x)$ on $W$ by $S(x)_W$. There exists a unique structure of truncated $(\mu , \nu )$-deformed $\phi$-coordinated $\Vcgl$-module on  $W$ such that the module map satisfies
$$
Y_{W}(T_{[n]}^+(u_1,\ldots ,u_n)\vac,z) =  S_{[n]}(ze^{u_1},\ldots , ze^{u_n})_W,
$$
where, in  parallel with \eqref{snixovi}, we write
$$
S_{[n]}(ze^{u_1},\ldots , ze^{u_n})_W=
  \prod_{i=1,\dots,n}^{\longrightarrow}
 S_i(ze^{u_i})_W\ts \R_{i\ts i+1}(1/z^2 e^{u_i+u_{i+1}})^{t_i}\ldots  \R_{i  n}(1/z^2 e^{u_i+u_{n}})^{t_i}.
$$
\end{kor} 

We have the following converse of    Corollary \ref{coro_conv_yo}.

\begin{kor}
Let $W$ be a truncated $(\mu , \nu )$-deformed $\phi$-coordinated $\Vcgl$-module such that
\beq\label{cond_def_rels1}
Y_W(t_{ij}^{(-1)}\vac,z)\in\om(W,z^{-1}W[[z^{-1}]])\quad\text{for all }1\leqslant i\leqslant j\leqslant N.
\eeq
There exists a unique structure of $\Yo$-module on $W$ such that
\beq\label{cond_def_rels2}
S(z)_W=Y_W(T^+(0)\vac,z).
\eeq
\end{kor} 

\begin{prf}
The $\mu $-commutativity \eqref{ficommprop} for the braiding map $\mu $  defined by \eqref{muh_formula}  is, essentially, of the same form as the reflection equation \eqref{refl_mult}, i.e., the defining relation for the $h$-twisted Yangian $\Yo$.
This observation, along with the fact that the condition \eqref{cond_def_rels1} ensures that the remaining defining relations  \eqref{refl_mult_2nd} for the $h$-twisted Yangian are satisfied by the right-hand side of \eqref{cond_def_rels2},  implies the corollary.
\end{prf}

The next corollary presents a  well-known result on the centrality of the coefficients of the {\em Sklyanin determinant}; see \cite[Cor. 4.3]{MRS}.
It is  an immediate consequence of Proposition \ref{mufixes} and the form of the $\mu$-commutativity \eqref{ficommprop} for the module map $Y_{\Yo}(\cdot,z)$ from Theorem \ref{thm_hyang}. Indeed, it suffices to observe that  \eqref{yyoformula} implies the identity
$$
Y_{\Yo}(\qdet T^+ (0)\vac,z)=\tr_{1,\ldots ,N}A^{(N)} S_{[N]}(z,ze^{-h},\ldots ,ze^{-(N-1)h})\ts D_1\ldots D_N .
$$

\begin{kor}
The coefficients of the series 
$$
\tr_{1,\ldots ,N}A^{(N)} S_{[N]}(z,ze^{-h},\ldots ,ze^{-(N-1)h})\ts D_1\ldots D_N \in\Yo[[z^{-1}]]
$$
belong to the center of the algebra $\Yo$.
\end{kor}

\section{Restricted modules for the generalized twisted \texorpdfstring{$h$}{h}-Yangian as   deformed \texorpdfstring{$\phi$}{phi}-coordinated \texorpdfstring{$\Vcgl$}{Vcrit(glN)}-modules}\label{section_04}

In this section, we consider the nontruncated version of     Definition \ref{defn_deformed_top},  so that the module map $Y_W(\cdot,z)$ is   allowed to have infinitely many positive and negative powers of the variable $z$. However, as with the definition of quantum vertex algebra \cite[Sect. 1.4.1]{EK5}, the series $Y_W(v,z)w$, when regarded modulo $h^n$, is required to possess a pole at $z=0$ for all $v,w$ and $n$. In parallel with Definition \ref{defn_deformed_top}, we   adapt the setting so that it features the  
associate $\phi(z_2,z_0)=z_2 e^{z_0}$.
To present the definition, we   need to introduce some notation.
Let $W$ be a topologically free $\CC[[h]]$-module and $ n $ a positive integer. Suppose that some series  $A=A(z_1,z_2)\in\om(W,W[[z_1^{\pm 1},z_2^{\pm 1} ]])$ can be expressed as
$$
A=B+  h^n C  \quad  \text{ for some }\quad
B\in\om(W,W((z_1,z_2))_h),\,  C \in\om(W,W [[z_1^{\pm 1},z_2^{\pm 1}]]).  
$$
Clearly, this can be expressed equivalently, without using the term $C$, as
$$
A-B\in h^n \om(W,W [[z_1^{\pm 1},z_2^{\pm 1}]]).
$$
To indicate   that the series $A $ possesses such a decomposition, we   write
$$
A(z_1,z_2)\in\om(W,W((z_1,z_2)) )\mod  h^n.
$$
Note that the substitution
$
B(z_1,z_2)\big|_{z_1=z_2 e^{z_0}}\big.$
is well-defined, even though the  substitution  $\textstyle A(z_1,z_2)\big|_{z_1=z_2 e^{z_0}} \big.$ does not need to be. 
In what follows,  having that in mind, we shall write
$$
A(z_1,z_2)\big|_{z_1=z_2 e^{z_0}}^{\text{mod }h^n}\big.=B(z_1,z_2)\big|_{z_1=z_2 e^{z_0}}\big. .
$$

\begin{defn}\label{def_def_mod_phi}
Let $(V,Y,\vac,\Sc)$  be a    quantum vertex algebra,  $(\mu  ,\nu  )$     a
   compatible pair on $V$ and
  $W$   a topologically free $\CC[[h]]$-module equipped with  a $\CC[[h]]$-module map
\begin{align*}
Y_W(\cdot ,z) \colon V\ot W&\to W((z))_h, \\
v\ot w&\mapsto Y_W(z)(v\ot w)=Y_W(v,z)w=\sum_{r\in\ZZ}  v_{r-1} w \ts z^{-r }.\non
\end{align*}
A pair $(W,Y_W)$ is said to be a {\em   $(\mu  ,\nu  )$-deformed $\phi$-coordinated $V$-module} if the   map $Y_W(\cdot, z)$ satisfies the {\em vacuum property},
\begin{align}
&Y_W(\vac, z)w=w\quad\text{for all }w\in W,\non\\
\intertext{the  {\em weak  $\nu $-associativity}:
for any elements $u,v \in V$  and $n \in\mathbb{Z}_{> 0}$
there exists a nonzero polynomial $p(z_1,z_2)\in\CC[z_1,z_2]$
such that we have}
&p(z_1,z_2)\ts Y_W(u,z_1)\ts Y_W(v,z_2)\in \om(W,W((z_1,z_2))) \mod h^n,\label{nuassocc1}\\
&\left(p(z_1,z_2)\ts Y_W(u, z_1)  Y_W(v,z_2)  \right)\big|_{z_1=z_2 e^{z_0}}^{\text{mod } h^n} \big.\non\\
&\qquad\big. - p(z_2 e^{z_0},z_2)\ts Y_W\big(Y^{\nu }_+ (u, z_0,z_2)v ,z_2\big)\in h^n \om(W,W[[z_0^{\pm 1},z_2^{\pm 1}]])   ,\label{nuassocc2}
\intertext{and the  {\em $\mu $-locality}:
for any $u,v\in V$ and $n\in\mathbb{Z}_{> 0}$ there exists
  a nonzero polynomial $q(z_1,z_2)\in\CC[z_1,z_2]$ such that we have }
 &\big(q(z_1,z_2)\ts Y_W(z_1)\big(1\otimes Y_W(z_2)\big)\big( \mu_{+}(z_1 /z_2,z_2)(u\otimes v)\otimes w\big)\big. 
\non\\
 &\qquad\big. - q(z_1,z_2)\ts Y_W(v,z_2) Y_W(u, z_1)  w \big)  \,
\in\,  h^n W[[z_1^{\pm 1},z_2^{\pm 1}]] \quad\text{for all }w\in W.\label{eslokaliti_hadica}
\end{align}
\end{defn}

\begin{rem}\label{rem_sing_2}
Note that the first term in the $\mu $-locality property  \eqref{eslokaliti_hadica} is well-defined due to the singularity constraint \eqref{sing_muh}. Moreover, the substitution $z_1=z_2e^{z_0}$ in the first term of the weak  $\nu $-associativity \eqref{nuassocc2} is well-defined modulo $h^n$, due to \eqref{nuassocc1}. Such a form of the weak associativity for the $\phi$-coordinated modules goes back to Li \cite[Def. 3.4]{Liphi}; see also \cite[Rem. 3.2]{Liphi} and \cite[Lemma 2.9]{LTW}.
\end{rem}

Motivated by the definition of the orthogonal twisted $h$-Yangian from \cite{MRS}, which we recalled in Section \ref{section_03}, we introduce the {\em generalized twisted $h$-Yangian} $\Xo$. It is defined as the $h$-adically completed algebra over the ring  $\CC[[h]]$ generated by the elements $b_{ij}^{(r)}$, where $i,j=1,\ldots ,N$ and $r\in\ZZ$, subject to  the defining relations 
\begin{align}
&\R(x/y)\ts B_1(x) \ts  \R(1/xy)^{t_1}\ts B_2(y)
=B_2(y)\ts  \R(1/xy)^{t_1}\ts B_1(x) \ts \R(x/y).\label{refl_mult_ext}
\end{align}
The matrix $B(u)$ is defined by
$$
B(u)=\sum_{i,j=1}^N e_{ij}\ot b_{ij}(u),\quad\text{where}\quad
b_{ij}(u)=\delta_{ij}+h\sum_{r\in\ZZ} b_{ij}^{(r)}u^{-r}.
$$ 

The identities \eqref{rmatrixprod3refl} and \eqref{rmatrixprod4}    extend analogously  to this setting, so that we have
\begin{align}
&\R_{nm}^{12}(x/y)\ts B^{13}_{[n]}(x) \ts  \RR_{nm}^{12\ts t}(1/xy) \ts B_{[m]}^{23}(y)
=\ts B_{[m]}^{23}(y)\ts  \RR_{nm}^{12\ts t}(1/xy) \ts B^{13}_{[n]}(x) \ts  \R_{nm}^{12}(x/y),\label{refl_general}\\
&B_{[n+m]}(x,y)=B^{13}_{[n]}(x) \ts  \RR_{nm}^{12\ts t}(1/xy) \ts B_{[m]}^{23}(y),\label{refl_general2}
\end{align}
where, for $x=(x_1,\ldots ,x_n)$, we write
\beq\label{bnixovi}
B_{[n]}(x)=
  \prod_{i=1,\dots,n}^{\longrightarrow}
 B_i(x_i)\ts \R_{i\ts i+1}(1/x_i x_{i+1})^{t_i}\ldots  \R_{i  n}(1/x_i x_{n})^{t_i}. 
\eeq

To connect  the generalized $h$-Yangian with the quantum affine   vertex algebra $\Vcgl$, we shall need the notion of restricted module. A $\Xo$-module $W$ is said to be {\em restricted} if it is a topologically free   $\CC[[h]]$-module such that the action $B_{[n]}(x)_W$ of the series \eqref{bnixovi} on $W$ satisfies
\beq\label{restr_def}
B_{[n]}(x )_W\in(\ndo\CC^N)^{\ot n}\ot\om (W,W((x_1,\ldots ,x_n))_h)\quad\text{for all }n=1,2,\ldots .
\eeq
 
\begin{rem}
In this remark, we provide   motivation for the above definition of restricted module.
First of all, we recollect  that, in \cite[Sect. 3.3]{MRS},  the authors discuss a certain matrix $\wvr{S}(u)$, which gives rise to an alternative definition of the orthogonal twisted $q$-Yangian. The  additive analogue of  $\wvr{S}(u)$ in the setting of this paper takes the form $\wvr{B}(u)\coloneqq T^+(u)T^*(-u)^t$, where the action of $T^*(u)$ is given by \eqref{tstar_formulaa} with $c=0$. Its powers are found by
$$
\wvr{B}_{[n]}(x)=
  \prod_{i=1,\dots,n}^{\longrightarrow}
 \wvr{B}_i(x_i)\ts R_{i\ts i+1}(e^{-x_i- x_{i+1}})^{t_i}\ldots  R_{i  n}(e^{-x_i- x_{n}})^{t_i},
$$
where $x=(x_1,\ldots ,x_n)$ and $n=1,2,\ldots ;$ cf. \eqref{bnixovi}. They exhibit the property
$$
\wvr{B}_{[n]}(x)=
T_{[n]}^+(x)\ts
\RRR_{[n]}^{ t}(e^{- x})\ts
T_{[n]}^*(-x)^{t_1,\ldots ,t_n} \quad\text{for}\quad
\RRR_{[n]}^{ t}(e^{-x})=\prod_{1\leqslant i<j\leqslant n} R_{ij}(e^{-x_i-x_j})^{t_i}, 
$$
where the product is taken in the lexicographical order on the set of pairs $(i,j)$.
Due to this property and the form of the $R$-matrix   \eqref{rplusg},    for all $n=1,2,\ldots$ we have
$$
\wvr{B}_{[n]}(x_1,\ldots ,x_n)_{\Ur}\in(\ndo\CC^N)^{\ot n}\ot\om (\Ur,\Ur((x_1,\ldots ,x_n))_h),
$$
which   motivated the requirement imposed by \eqref{restr_def}.
\end{rem}

As before, let  $(\mu , \nu )$ be the compatible pair from Proposition \ref{nuhatprop}.
\begin{thm}\label{thm_hyang_ext}
Let $W$ be a restricted  $\Xo$-module.
There exists a unique structure of   $(\mu , \nu )$-deformed $\phi$-coordinated $\Vcgl$-module on   $W$   such that the module map satisfies
\beq\label{modulemapw}
Y_{W}(T_{[n]}^+(u_1,\ldots ,u_n)\vac,z) =  B_{[n]}(x_1,\ldots , x_n)_W \big|_{x_1=z e^{u_1},\ldots ,\ts x_n=ze^{u_n}}.
\eeq
\end{thm}

\begin{prf}
Suppose $W$ is a restricted  $\Xo$-module.
First of all, we observe that the substitution on the right-hand side of \eqref{modulemapw} is well-defined due to \eqref{restr_def}. In addition, the coefficients of all $u_1^{a_1}\ldots u_n^{a_n}$ of the matrix entries in the resulting expression belong  to $\om(W,W((z))_h)$. Thus, we conclude that $Y_W(v,z)w$ belongs to   $W((z))_h$ for any $v\in\Vcgl$ and $w\in W$, as required. 
To simplify the notation, we shall  denote the right-hand side of  \eqref{modulemapw} more briefly by $B_{[n]}(z e^{u_1},\ldots , z e^{u_n})_W =B_{[n]}(z e^{u } )_W$.

Next, we show that \eqref{modulemapw}, together with $Y_W(\vac,z)=1_W$, defines a $\CC[[h]]$-module map. Let $i$, $n$ be positive integers such that $i<n$. 
Introduce the index sets  
$$I_1=\left\{1,\ldots ,i-1\right\} ,\quad I_2=\left\{i,i+1\right\} ,\quad I_3=\left\{i+2,\ldots ,n\right\},\quad I_4=\left\{n+1\right\} .$$
We shall  indicate the corresponding tensor factors of
$$
\overbrace{(\ndo\CC^N)^{\ot (i-1)}}^{1}\ot
\overbrace{(\ndo\CC^N)^{\ot 2}}^{2}\ot
\overbrace{(\ndo\CC^N)^{\ot (n-i-1)}}^{3}\ot
\overbrace{W}^{4}.
$$ 
by writing  the numbers $1$--$4$ in the superscripts.
By using the decomposition in \eqref{refl_general2}, we can express $B_{[n]}(x )_W = B_{[n]}(x_1,\ldots ,x_n)_W$ as
$$
B_{[i-1]}^{14}(x_1,\ldots,x_{i-1})_W\ts
\RR^{ 12\ts t}_{i-1\ts 2}\ts
B_{[2]}^{24}(x_i,x_{i+1})_W\ts
\RR^{ 13\ts t}_{i-1\ts n-i-1}\ts
\RR^{23\ts t}_{ 2\ts n-i-1}\ts
B_{[n-i-1]}^{34}(x_{i+2},\ldots,x_{n})_W,
$$
where $x=(x_1,\ldots ,x_n)=(y_1e^{u_1},\ldots ,y_ne^{u_n})$  is the $n$-tuple of variables  and
$$
\RR_{|I_a|\ts |I_b|}^{ ab\ts t}  = \prod_{i\in I_a}^{\longrightarrow}
\prod_{j\in I_b}^{\longrightarrow}   \R_{ij}(1/x_i x_j)^{t_i}.
$$
Hence, by using the Yang--Baxter equation \eqref{YBE} and the reflection equation \eqref{refl_mult_ext}, we get
\begin{align}
\R_{i\ts i+1}\ts B_{[n]}(x )_W
=&\, B_{[i-1]}^{14}(x_1,\ldots,x_{i-1})_W\ts
P_{i\ts i+1}\ts \RR^{ 12\ts t\ts \prime}_{i-1\ts 2} P_{i\ts i+1}\ts
\ts \cev{B}_{[2]}^{24}(x_{i },x_{i+1})_W\non  \\
&\times\RR^{ 13\ts  t}_{i-1\ts n-i-1}\ts P_{i\ts i+1}\ts
 \RR^{23\ts t\ts \prime}_{ 2\ts n-i-1}\ts P_{i\ts i+1}\ts
B_{[n-i-1]}^{34}(x_{i+2},\ldots,x_{n})_W\ts \R_{i\ts i+1},\label{tmp1}
\end{align}
where $\R_{i\ts i+1}=\R_{i\ts i+1}(x_i /x_{i+1})$, the prime in $\RR_{|a|\ts |b|}^{ ab\ts t\ts \prime}$ indicates that the variables $x_i$ and $x_{i+1}$ are swapped,  $P_{i\ts i+1}$ stands for the action of the permutation operator
\eqref{permutation_operator}
on the tensor factors $i$ and $i+1$ and we use the left arrow on the top of the symbol to indicate that its factors come in the reversed order:
$$
\cev{B}_{[2]}^{24}(x_{i },x_{i+1})_W=B_{i+1\ts n+1}(x_{i+1})_W\ts \R_{i\ts i+1}(1/x_i x_{i+1})^{t_i}\ts B_{i \ts n+1}(x_{i })_W.
$$
 It remains to observe that by applying the substitutions $y_1=z,\ldots ,y_n=z $ to the left-hand side (resp. right-hand side) of \eqref{tmp1}, we obtain the image of $\R_{i\ts i+1}(e^{u_i -u_{i+1}})T_{[n]}^+(u)\vac$ (resp. $P_{i\ts i+1} T_{[n]}^{+\ts \prime}(u)P_{i\ts i+1} \R_{i\ts i+1}(e^{u_i -u_{i+1}})\vac$) under \eqref{modulemapw}, where prime in $T_{[n]}^{+\ts \prime}(u)$   means that   $u_i$ and $u_{i+1}$ are swapped. Finally, we observe that the defining relations \eqref{rtt} for $\Ur=\Vcgl$ can be written equivalently in terms of the $R$-matrix $\R(e^u)$ as
$$
\R(e^{u-v})\ts T_{1}^+(u)\ts T_2^+ (v)=  T_2^+ (v)\ts T_{1}^+(u)\ts \R(e^{u-v}).
$$
Thus, our discussion shows that the ideal of   defining relations \eqref{rtt}   is in the kernel of \eqref{modulemapw}, so  that the $\CC[[h]]$-module map $Y_W(\cdot, z)$ is well-defined.

Let us prove that  $Y_W(\cdot ,z)$ satisfies the weak $\nu$-associativity property  \eqref{nuassocc1}.   Set
$$ 
x=(z_1e^{u_1},\ldots,z_1 e^{u_n}),\quad y=(z_2 e^{v_1},\ldots ,z_2 e^{v_m}),\quad
u=(u_1,\ldots ,u_n),\quad v=(v_1,\ldots ,v_m).
$$ 
Due to \eqref{csym_cons}, Equality \eqref{refl_general2}
implies
\beq\label{tmpp2}
B^{13}_{[n]}(x)_W   \ts B_{[m]}^{23}(y)_W=\iotadva\RR_{nm}^{12\ts *}(1/xy)  \cdotrl B_{[n+m]}(x,y)_W.
 \eeq
As $W$ is a restricted module, the coefficients of all $u_1^{a_1}\ldots u_n^{a_n}v_1^{b_1}\ldots v_m^{b_m}$ in the matrix entries of  $B_{[n+m]}(x,y)_W$ belong to $\om(W,W((z_1,z_2))_h)$; recall \eqref{restr_def}. Consider the other term on the right-hand side, $\iotadva\RR_{nm}^{12\ts *}(1/xy) $.  The   property \eqref{unitrig} implies  
\beq\label{r-invers}
\R_{12} (z)^{-1}=-\frac{ze^h}{(1-ze^h)(1-ze^{-h})} \R_{21} (1/z).
\eeq
Therefore, it is evident from \eqref{exp} that for any integers $c_1,c_2\geqslant 0$ and  $k> 0$ there exists a nonzero polynomial $q(z)$ such that the coefficients of $w_1^{c_1}w_2^{c_2}$ in the matrix entries of
$  q(z)  \iotaze \R_{12} (ze^{w_1+w_2})^{-1}$ 
belong to $\CC[z ,h]$ modulo $h^k$. By this observation, one deduces that for any integers   $a_1,\ldots ,a_n,b_1,\ldots ,b_m\geqslant 0$ and $k>0$, there exists a nonzero polynomial $p(z_1,z_2)$ such that
the coefficients of   all
\beq\label{varrs}
 u_1^{a_1^\prime}\ldots u_n^{a_n^\prime}v_1^{b_1^\prime}\ldots v_m^{b_m^\prime}  , \text{ where } 0\leqslant a_i^\prime\leqslant a_i,\, 
 0\leqslant b_j^\prime\leqslant b_j,\ i=1,\ldots ,n,\,j=1,\ldots ,m ,
\eeq
in the matrix entries of
$ 
p(z_1,z_2)\iotadva\RR_{nm}^{12\ts *}(1/xy) 
$ 
belong to $\CC [z_1,z_2 ,h]$ modulo $h^k$. The preceding discussion, together with     \eqref{tmpp2}, implies that
the coefficients of   \eqref{varrs}
in the matrix entries of
\beq\label{choosepoly}
p(z_1,z_2)\ts B^{13}_{[n]}(x)_W   \ts B_{[m]}^{23}(y)_W
=p(z_1,z_2)\ts Y_W(T_{[n]}^+(u )\vac,z_1)\ts Y_W(T_{[m]}^+(v )\vac,z_2)
\eeq
belong to 
$\om(W,W((z_1,z_2)) )$ modulo $h^k$, so we conclude that \eqref{nuassocc1} holds.

To finish the proof of  the weak $\nu$-associativity, it remains to check the requirement imposed by \eqref{nuassocc2}. Replacing the vector $u$ (resp. $v$) in the second summand of \eqref{nuassocc2} by the series $T_{[n]}^+(u_1,\ldots ,u_n)$ (resp. $T_{[m]}^+(v_1,\ldots ,v_m)$), then using the explicit expression \eqref{ynu_expl} for the action of $Y^\nu_+(\cdot,z_0,z_2)$ and, after all, by employing \eqref{modulemapw}, we obtain
\begin{align}
&p(z_2 e^{z_0},z_2)\ts Y_W\big(Y^{\nu }_+ (T_{[n]}^+(u)\vac , z_0,z_2)T_{[m]}^+(v ) \vac,z_2\big)\label{tmpxpr1n}\\
=&\,p(z_2e^{z_0},z_2)\iotaopdvanula \RR_{nm}^{12\ts *}(1/z_2^2 e^{z_0+u+v}) \cdotrl B_{[n+m]}(z_2e^{z_0+u_1},\ldots ,z_2e^{z_0+u_n},z_2e^{v_1},\ldots ,z_2e^{v_m})_W,\non
\end{align}
 where $p$ is the same polynomial as in \eqref{choosepoly}.
On the other hand, by the earlier discussion and, in particular, \eqref{tmpp2}, the coefficients of   \eqref{varrs}
in the matrix entries of \eqref{choosepoly}, when regarded modulo $h^k$, coincide with the same coefficients in
$$\left(\iotadva p(z_1,z_2)\RR_{nm}^{12\ts *}(1/xy)  \right)\cdotrl B_{[n+m]}(x,y)_W.$$
Therefore,  by \eqref{choosepoly}, the coefficients of  \eqref{varrs}
in the matrix entries of
\begin{align}
& \left(\left(\iotadva p(z_1,z_2)\RR_{nm}^{12\ts *}(1/xy)  \right)\cdotrl B_{[n+m]}(x,y)_W\right)\Big|_{z_1=z_2e^{z_0}}^{\text{mod }h^k}\non\\
=&\left( \iotadva p(z_1,z_2)\RR_{nm}^{12\ts *}(1/xy) \right) \Big|_{z_1=z_2e^{z_0}}^{\text{mod }h^k} \cdotrl B_{[n+m]}(x,y)_W \Big|_{z_1=z_2e^{z_0}}^{\text{mod }h^k}\label{tmpxpr2n}
\end{align}
and
\beq\label{tmpxpr3n}
 \left(p(z_1,z_2)\ts Y_W(T_{[n]}^+(u )\vac,z_1)\ts Y_W(T_{[m]}^+(v )\vac,z_2)\right)\Big|_{z_1=z_2e^{z_0}}^{\text{mod }h^k} 
\eeq
coincide. Finally, by comparing \eqref{tmpxpr1n} and \eqref{tmpxpr2n}, we conclude that  the   coefficients of \eqref{varrs} in the matrix entries of  \eqref{tmpxpr1n} and \eqref{tmpxpr3n} coincide modulo $h^k$, which implies   the weak $\nu$-associativity property
\eqref{nuassocc2}.

Let us prove the $\mu$-locality \eqref{eslokaliti_hadica}. 
Replacing the vector $u$ (resp. $v$) in the first summand of \eqref{eslokaliti_hadica} by   $T_{[n]}^+(u_1,\ldots ,u_n)$ (resp. $T_{[m]}^+(v_1,\ldots ,v_m)$) and then using \eqref{muh_formula} and \eqref{modulemapw} we get
\begin{align*}
&\RR_{nm}^{12\ts *}(1/z_1 z_2 e^{ u+v}) 
\cdotlr\left(\R_{nm}^{12}(z_1 e^{u-v}/z_2)\ts B_{[n]}^{ 13}(z_1 e^u)_W\ts \right.\non\\
&\left.\times\RR_{nm}^{12\ts t}(1/z_1 z_2 e^{ u+v}) \ts B_{[m]}^{ 23}(z_2 e^v)_W\ts \R_{nm}^{12}(z_1 e^{u-v}/z_2)^{-1}\right),
\end{align*}
where the $R$-matrices are expanded so that the coefficients of  all $u_1^{a_1}\ldots u_n^{a_n}v_1^{b_1}\ldots v_m^{b_m}$ in their matrix entries belong to
$ \CC[z_1^{\pm 1}]((z_2))[[h]]$;
recall \eqref{sing_muh}. By using the generalized reflection equation \eqref{refl_general}, we rewrite the above expression as
\begin{align*}
&\RR_{nm}^{12\ts *}(1/z_1 z_2 e^{ u+v}) 
\cdotlr\left(\left(B_{[m]}^{ 23}(z_2 e^v)_W\ts
\RR_{nm}^{12\ts t}(1/z_1 z_2 e^{ u+v})  
  \right.\right.\non\\
&\left.\left.
\times B_{[n]}^{ 13}(z_1 e^u)_W\ts
\R_{nm}^{12}(z_1 e^{u-v}/z_2)\right)
 \R_{nm}^{12}(z_1 e^{u-v}/z_2)^{-1}\right).
\end{align*}
Due to the form of the $R$-matrix inverse in \eqref{r-invers}, we conclude that 
for any integers   $a_1,\ldots ,a_n,b_1,\ldots ,b_m\geqslant 0$ and $k>0$, there exist   nonzero polynomials $q_1(z_1,z_2)$ and $q_2(z_1,z_2)$ such that
the coefficients of    \eqref{varrs}
in the matrix entries of
$$
q_1(z_1,z_2)\ts \RR_{nm}^{12\ts *}(1/z_1 z_2 e^{ u+v}) 
\fand
q_2(z_1,z_2)\ts \R_{nm}^{12}(z_1 e^{u-v}/z_2)^{-1},
$$
when regarded modulo $h^k$,
lie in $ \CC[z_2^{\pm 1},h]((z_1)) $. Thus, with $q (z_1,z_2)=q_1(z_1,z_2)q_2(z_1,z_2)$, this implies that the coefficients of \eqref{varrs} in the matrix entries of
\begin{align*}
&q(z_1,z_2)\ts \RR_{nm}^{12\ts *}(1/z_1 z_2 e^{ u+v}) 
\cdotlr\left(\left(B_{[m]}^{ 23}(z_2 e^v)_W\ts
\RR_{nm}^{12\ts t}(1/z_1 z_2 e^{ u+v})  
  \right.\right.\non\\
&\left.\left.
\times B_{[n]}^{ 13}(z_1 e^u)_W\ts
\R_{nm}^{12}(z_1 e^{u-v}/z_2)\right)
 \R_{nm}^{12}(z_1 e^{u-v}/z_2)^{-1}\right) 
\end{align*}
coincide with the coefficients of \eqref{varrs} in the matrix entries of
\beq\label{theothersummand}
 q(z_1,z_2)  \ts B_{[m]}^{ 23}(z_2 e^v)_W\ts
  B_{[n]}^{ 13}(z_1 e^u)_W 
\eeq
modulo $h^k$. To finish the proof, it suffices to observe that \eqref{theothersummand}
is equal to 
$$
q(z_1,z_2)  \ts Y_W(T_{[m]}^{+ 23}(v)\vac,z_2)\ts
  Y_W(T_{[n]}^{+ 13}(u)\vac,z_1), 
$$
i.e., in other words, it
 corresponds to the second summand in the $\mu$-locality property
\eqref{eslokaliti_hadica}. Hence, we conclude that the map $Y_W(\cdot,z)$ satisfies the $\mu$-locality, as required.
\end{prf}

 \section*{Acknowledgment}
L.B. is member of Gruppo Nazionale per le Strutture Algebriche, Geometriche e le loro Applicazioni  (GNSAGA) of the Istituto Nazionale di Alta Matematica (INdAM).
This work has been supported  by Croatian Science Foundation under the project UIP-2019-04-8488. Furthermore, this work was supported by the project ``Implementation of cutting-edge research and its application as part of the Scientific Center of Excellence for Quantum and Complex Systems, and Representations of Lie Algebras'', PK.1.1.02, European Union, European Regional Development Fund.

\end{document}